\documentclass[12pt]{article}
\usepackage{amsmath,amsfonts,amssymb,amsthm,amscd}
\usepackage{a4wide}
\title{{\it Galois descent in  Galois theories}}
\author{Daniel Bertrand\thanks {~Institut de Math\'ematiques de Jussieu. - Adresse, mots-clefs et classification en fin de texte.}}
\newtheorem{Theorem}{Theorem}[section]

\newcommand{\C}{\mathbb C}
\newcommand{\R}{\mathbb R}
\newcommand{\Q}{\mathbb Q}
\newcommand{\Z}{\mathbb Z}
\newcommand{\F}{\mathbb F}

\newcommand{\K}{\mathbb K}
\newcommand{\G}{\mathbb G}

\date{March 21,  2010  (revised May 9)}
\begin{document}
\maketitle

\noindent
{\it Abstract} : inspired by Kummer theory on abelian varieties, we give similar looking descriptions of the Galois groups occuring in the differential Galois theories of Picard-Vessiot, Kolchin and Pillay, and mention some arithmetic applications.

\medskip
\noindent
{\it R\'esum\'e} : guid\'es par la th\'eorie de Kummer sur les vari\'et\'es ab\'eliennes et motiv\'es par quelques applications arithm\'etiques, nous donnons des descriptions d'apparences similaires des groupes de Galois issus des th\'eories de Galois diff\'erentielles de Picard-Vessiot, Kolchin et Pillay.

\bigskip

\centerline{***}

\bigskip

The topic I had been assigned by the organizers of the Luminy September 09 School was  ``Algebraic $D$-groups and non-linear differential Galois theories". The present account is written in an applied maths spirit : how to compute the Galois groups, and what for\;? Thus, we start with a motivating  question which, in accordance with the theme of the School, comes from diophantine geometry. We then describe the Galois groups of the various theories under study, in terms that bear a strong similarity. Finally, we apply this description to the study of exponentials and logarithms on abelian schemes. 

\medskip
A general argument of Galois descent occurs along the text, hence the title\footnote{~also borrowed from a set of talks at the Durham July 09  Conference on model theory. I thank the organizers of both Luminy and Durham meetings for having offered me these opportunities to develop this point of view.} of these notes; its number theoretic prototype, given by Kummer theory, is recalled in an Appendix to the paper. 

\medskip
Although the presentation is sometimes novel, the results described here are not new. For original sources, we refer the reader to \cite{vdP-S} for the Picard-Vessiot theory, \cite{Pi} for Kolchin's and Pillay's theories,   and  to \cite{An} and \cite{B-P} for the applications to algebraic independence. Actually, this text may serve as an introduction to the survey \cite{BM}, which is itself an introduction to the latter papers (and to the descent argument in the non-linear case). 

\section{Ax-Schanuel}

\subsection{The multiplicative case}

The well-known Schanuel conjecture asserts that if $x = \{x_1, ..., x_n\}$ is a ``non-degenerate"  $n$-tuple of complex numbers whose image under the standard exponential function $exp$ is denoted by $y = \{y_1 = exp(x_1), ..., y_n = exp(x_n)\}$, then $tr.deg._{\Q} \Q(x,y) \geq n$. The expression {\it non-degenerate} will occur under several acceptions in these notes. Here, it means that 
$$\forall (m_1, .., m_n) \in \Z^n \setminus 0,  m_1x_1 + ... + m_n x_n \neq 0,$$
or equivalently, that for any proper algebraic subgroup $H$ of the algebraic torus  $G = \G_m^n$, the complex point $x$ of the Lie algebra $LG $ of $G$ does not lie in the Lie algebra $LH$ of $H$. 

\medskip
In 1970, Ax \cite{Ax1} proved a functional version of the conjecture, which, in an analytic setting, may be phrased as follows. Let
$$x = (x_1,..., x_n) \in (\C\{\{z_1, .., z_t\}\})^n$$
be a $n$-tuple of convergent power series in $t$ variables. For each $i$, $ y_i(z) := exp(x_i)$ lies in $ \C\{\{z_1, .., z_t\}\}^*$, and we set  $  y = exp(x) \in  (\C\{\{z_1, .., z_t\}\}^*)^n$.
Assume that $x$ is {\it non-degenerate}, in the sense that
$$\forall  (m_1, .., m_n)  \in \Z^n \setminus 0,  m_1x_1 + ... + m_n x_n \notin \C,$$
or equivalently, that  for any proper algebraic subgroup $H$ of the algebraic torus  $G = \G_m^n$ and any constant point $\xi \in LG(\C)$, 
$$ x  \notin \xi +  LH.$$
Then,   $tr.deg. (\C(x,y)/\C) \geq n +  \mu$, where $\mu$ denotes the rank of the functional jacobian ${Dx \over Dz} \in Mat_{t,n}(\C\{\{z_1, .., z_t\}\})$ of $x$.

\medskip

Let  $K \simeq \C(z_1, ..., z_\mu)$ be the field generated by the principal variables. In order to check the above lower bound, it suffices to show that $tr.deg. (K(x,y) /K )\geq n$ (and the two statements are actually equivalent). Furthermore, choosing a sufficiently general line in $\C^\mu$,  it suffices to check the latter inequality   when $\mu = 1$. Using the differential equation satisfied by $exp$, we can therefore view Ax's theorem as a corollary of the following differential algebraic statement.   Let $K = \C(z)^{alg}$ be the algebraic closure of $\C(z)$, endowed with the (unique) extension $\partial$ of the derivation $ \frac{d}{dz}$, and let $ ({\cal K}, \partial)$ be a differential extension of $(K, \partial)$, with same constant field ${\cal K}^\partial = \C$. Let further
$(x, y) \in ({\cal K}Ê\times {\cal K}^*)^n$ satisfy : $Ê \partial y /y = \partial x$ (where derivations are taken coordinate-wise). Assume that for any proper algebraic subgroup $H$ of $G = \G_m^n$,  $x$ does not lie in $LH  + LG(\C)$.   Then 
$deg.tr. K(x,y) /K \geq n$.

\subsection{The constant case}
Two years later, Ax \cite{Ax2}  extended his results to more general algebraic groups (see also \cite{Ki}). For instance, making the same analysis as above, we may rephrase Theorem 3 of his paper (actually written in a formal setting and under a slightly stronger hypothesis on $G$) as follows. 

\medskip
Let    $K = \C(z)^{alg},  \partial = \frac{d}{dz}$ and  ${\cal K}$ be as above, and  let 
$G$ be a commutative algebraic group {\it defined over $\C$},  with no additive quotient. In other words, $G$ is a semi-abelian variety defined over $\C$, or more generally, a quotient of its universal vectorial extension. The Lie algebra $LG$ of $G$ is a  vector space over $\C$, so that  there is a unique differential operator  $\nabla_{LG}$ on $LG({\cal K})$, whose solution space is $LG(\C)$ (to define $\nabla_{LG}(x)$, choose any basis of $LG$ {\it over $\C$}, and take the $\partial$-derivatives  of  the coordinates of $x$; the outcome is  independent of the chosen basis). The  exponential map $exp_G : LG(\C) \rightarrow G(\C)$ is a morphism of commutative  Lie groups, admitting as kernel  a discrete subgroup $\Omega_G$ of $LG(\C)$, and one can consider its inverse $\ell n_G$ as a multivalued function.  For any analytic function $x(z)$ with values in $LG(\C)$, $y(z) := exp_G(x(z))$ is a well defined analytic function with values in $G(\C)$. For any analytic function $y(z)$ with values in $G(\C)$, $\nabla_{LG}Ê\circ \ell n_G(y)$ is also well-defined, since $\Omega_G$ is killed by $\nabla_{LG}$. Its explicit expression enables us to extend $\nabla_{LG} \circ  \ell n_G $ to a group homomorphism $\partial \ell n_G : G({\cal K}) \rightarrow LG({\cal K})$.  This is the {\it logarithmic derivative} of $G/\C$, which we describe in a more algebraic way in \S 2, then for non constant groups  in \S 3 - and again  in the above style in \S 4.1. Notice that when $x$ and $y$ have an analytic meaning, the relation $\partial \ell n_G(y) = \nabla_{LG}(x)$ is equivalent to the existence of a point $\xi \in LG(\C)$ such that $y  = exp_G(x -\xi)$. 

\medskip

Exactly as in \S 1.1, Ax's theorem then reads as follows : let $(x,y) \in (LG \times G)({\cal K})$ satisfy $\partial \ell n_G (y) = \nabla_{LG}(x)$, and suppose that  $x$ is non-degenerate : for any   proper
algebraic subgroup $H/\C$
of $G$,  $x \notin LH + LG(\C)$. 
Then, $ tr.deg. (K(x,y)/K) \geq dim  G$.
\bigskip

In these notes, we will show that   differential Galois theories  provide proofs of Ax's theorem under the following {\it restrictions} :

$(L)$ either  $y \in G(K)$ (in which case we can apply the linear Picard-Vessiot theory);

$(E)$ or $x \in LG(K)$ (in which case we can apply the non-linear theory of Kolchin).

\noindent
But the interesting point about  these Galois approaches is that in  fact, they then provide an {\it extension} of these results to the case of {\it non-constant} algebraic groups, where  $G$ will only be defined over $K$. In the second situation, this is made possible by Pillay's generalization  of Kolchin's theory (although the initial proof given in \cite{B-P} uses a different method). See Theorems 4.1 to 4.4 for the outcome in the case of abelian varieties.

\subsection{Motivations}

The Manin-Mumford conjecture was proved by Raynaud in 1984, and has known since then a remarkable number of interesting new proofs. Based on work of Bombieri, Pila, Wilkie, and Zannier, Pila recently obtained another one \cite{Pa}, where the strategy of \cite{PZ} is  combined with  Ax's theorem\footnote{~That Ax's earlier version on tori  \cite{Ax1}  can play a similar role for the multiplicative analogue of Manin-Mumford, had already  been observed in \cite{PZ}, Final Remark 2.} on abelian varieties over $\C$. By a general argument (see \cite {Bog}, Thm. 1 and proof of Thm. 4; also \cite{PZ}), the conjecture reduces to :

\medskip
\noindent
{\bf Manin-Mumford} (key point): {\it Let $A/\Q^{alg}$ be an abelian variety. An algebraic subvariety $X/\Q^{alg}$  of $A$ passes through finitely many torsion points of $A$, unless $X$  contains a translate of a non-zero abelian subvariety of $A$. }

\medskip

We now sketch Pila's approach : first, as in \cite{PZ}, write $A(\C) = LA(\C)/\Omega_A \simeq \R^{2g}/ \Z^{2g}$, so that the torsion points become the rational points of the box  $[0,1[^{2g}$  while $X$ pulls back to a complex analytic subvariety $\cal X$ of $LA$. By $o$-minimality, $\cal X$ meets  $<< T^{\epsilon}$ rational points of denominator $\leq T$,  {\it outside of the real semi-algebraic subvarieties ${\cal W}$ of positive dimension  it contains}.  But back to $A(\Q^{alg})$, any such torsion point generates many others\footnote{~i.e. more than $T^\delta$, where $\delta > \epsilon$. The  Kummer theory described in the Appendix  would similarly yield large Galois orbits  for the division points in the Mordell-Lang conjecture.} by Galois action,  so their orders are bounded.  
\medskip

To conclude, we must control  the possible irreducible  complex algebraic subvarieties   $W$ of positive dimension in ${\cal X}$. Assuming that $X$ contains no translate of a non zero abelian subvariety of $A$, we claim that no such   $W$ exists. Assuming otherwise, consider the function field $K = \C({W})$, and let $x \in LA(K)$ be a generic point over $\C$ of   $W$.  Since $exp_A(W) \subset X$, the transcendence degree of $y = exp_A(x)$ over $\C$ is $< dim(A)$.   Ax's theorem on the constant abelian variety $G = A$ (in its original formulation, or in the above one, using A. Pillay's remark that it suffices to check the claim when $W$ is a curve)  then implies that   $x$ lies in  $\xi  + LA'$, for  some abelian subvariety  $A'$  of $A$, with $0 \neq A'  \subsetneq A$,  and some  $\xi  \in LA(\C)$. Set  $\eta = exp_A(\xi)$, and notice that $x' = x - \xi \in LA'(K)$ is still a generic point over $\C$ of the irreducible algebraic variety $W' = W - \xi$, which is therefore contained in $LA'$. The image of $W'$ under $exp_{A'}$ is contained in the  intersection $X'$ of $X - \eta$ and $A'$, which is a  subvariety of $A'$ containing no translate of a non zero abelian subvariety. Since the inverse image ${\cal X}' \subset LA'$ of $X'$ under $exp_{A'}$ contains $W'$,  the proof can now be concluded by induction on the dimension of $A$. 
 We point out  that Ax's theorem was here  used  only in the $(E)$ setting.

\bigskip
In his unpublished note \cite{Pk}, Pink extended the conjecture to a relative context, including the following case :

\medskip
\noindent
{\bf Relative Manin-Mumford} (over curves) : {\it let $X$ be the image of a non-torsion section of an abelian scheme ${\cal A}/S$ of relative dimension $\geq 2$ over a curve $S/\C$. Assume that $X$ is not contained in a translate of an elliptic subscheme of ${\cal A}/S$. Then, $X$  should meet finitely many of the torsion points of the various fibers of ${\cal A}/S$.}

\medskip
So, we here have an abelian variety $A$ over $K = \C(S)$. It need not come from $\C$, but one may hope that again,  an Ax-type theorem, now over a non-constant algebraic group, will help. And indeed, in their work on  the conjecture, Masser and Zannier  do appeal to such an algebraic independence statement, though now in the $(L)$ setting : see \cite{MZ}, p. 493, line 14, for  the test case of   the square of an elliptic scheme.

\section{ Picard-Vessiot \& Kolchin}

The differential Galois theories attached to these names (the second one generalizing the first one) concern constant algebraic groups. More precisely, let $(K, \partial)$ be a differential field, with an algebraically closed constant field   $K^\partial := C$ of characteristic $0$. We fix a differential closure $\hat K$ of $K$; in particular, $\hat K^\partial = C$. Let further
$G$ be a connected algebraic group, possibly non commutative, defined over $C$; the Picard-Vessiot theory  concerns affine algebraic groups  $G \subset GL_{n/C}$. We denote by $LG
$ the Lie algebra of $G$, and set : $G_K = G \otimes_C K$. The main point in the approach of Kolchin and his school (see in particular \cite{Ko}) is the existence of the  logarithmic derivative of $G$, a canonical differential algebraic map
 $$\partial \ell n_G  : G \rightarrow LG,$$
 which, at the level of $K$-rational points, can be described as follows.

\medskip
Any point $y \in G(K)$ provides a  derivation $\delta_Cy$ on the local ring ${\cal O}_{G/C, y}$, with values in $K$, via the formula $\delta_Cy(f_C) = \partial (f_C(y))$. By $K$-linearity, we can extend $\delta_Cy$  to a $K$-linear derivation $\partial y$ on  ${\cal O}_{G_K, y}$. In preparation for the next section,  we repeat the definition of $\partial y$ in the framework of \cite{Bu}, \cite{Pi}. First, $\partial$ extends uniquely to a derivation $D^0_\partial$ of the structure sheaf ${\cal O}_{G_K}$, killing ${\cal O}_{G/C}$. For instance,  for $G$ affine, we have $D^0_\partial = 1 \otimes \partial$ on $K[G] = C[G]Ê\otimes_C K$, while $\partial y  = \delta_C y \otimes 1$. 
 For any point $y \in G(K)$,  the formulae $D^0_{\partial, y}(f) = (D^0_\partial f)(y)$ and $(\delta y)(f) = \partial(f(y))$ define two derivations on the local ring ${\cal O}_{G_K, y}$, with values in $K$, which are only $C$-linear. But since both reduce to $\partial$ on $K$, their difference  $\delta y  - D^0_{\partial, y}$ vanishes on $K$, hence  is $K$-linear - and clearly coincides with $\partial y$ : $\partial(f(y)) =    \partial y(f)+ (D^0_\partial f)(y)$. 
 
 Now, such a  $K$-linear derivation  $\partial y : {\cal O}_{G_K, y} \rightarrow K$ can  be viewed as
  an element  of the tangent space $T_yG_K$ of $G_K$ at $y$. Pulling $\partial y$  back to the tangent space $LG_K$ of $G_K$ at the origin by the differential of right translation by $y^{-1}$, we obtain the {\it logarithmic derivative} $\partial \ell n_G(y) \in LG(K)$  of $y$ with respect to the standard extension $D^0_\partial$ of $\partial$.  In the affine case $G \subset GL_{n/C}$, this is given by
$$G(K) \ni y \mapsto \partial \ell n_G(y) = (\partial y) y^{-1} \in LG(K), $$
where $ \big((\partial y)_{ij}\big) = \big(\partial (y_{ij})\big) \in T_yG \subset T_y(Mat_{n,n}) \simeq Mat_{n,n}$. In particular, $\partial \ell n_{{\bf G}_m} y = {\partial y \over y} := \partial \ell n y ~,~  \partial \ell n_{{\bf G}_a}  y = \partial y.$ 
 
These formulae make sense over any differential extension $({\cal K}, \partial)$  of $(K, \partial)$. For any such ${\cal K}$, we write
$$G^\partial({\cal K}) = \{ y \in G({\cal K}), \partial \ell n_G(y) = 0\},$$
and point out  that in the present constant case, we have $G^\partial(\hat K) = G(C)$. We also note that $\partial \ell n_G$ is surjective at the level of $\hat K$-rational points, cf. \cite{Ko}, Prop. 6.

\medskip

Given $a \in LG(K)$, the Picard-Vessiot-Kolchin theory  studies the differential extension $K(y)/K$, where $y$ is a
 solution in $G(\hat K)$ of  
$$ \partial \ell n_G(y) = a,$$
and its Galois group
$$\Gamma_a = Aut_\partial(K(y)/K)    := \{\sigma \in Aut(K(y)/K),  \sigma \partial = \partial \sigma\}.$$  
The logarithmic derivative is a cocycle for the adjoint action of $G$ on $LG$:
$$\partial \ell n_G(uv) = \partial \ell n_G u + u (\partial \ell n_G v) u^{-1},$$
or equivalently : $\partial \ell n(u^{-1}v) = u^{-1}(-\partial \ell n u + \partial \ell n v) u$. Therefore,

 i) two solutions $y, \tilde y$ satisfy $y^{-1} \tilde y  \in G^\partial(\hat K) = G(C)$, so, the field $K(y)$ does not depend on the choice of $y$, and for a given $y$, the Galois group admits a natural embedding $\rho$ into $G(C)$ :  $\forall \sigma \in \Gamma_a$, $y^{-1}\sigma y =  \rho(\sigma) \in G^\partial(\hat K) = G(C)$. If we replace $y$ by another solution $yc, c \in G(C)$, then, $\rho$ is changed into $c^{-1} \rho c$.

 ii) Consider the ``affine" action of  $G(K)$ on $LG(K)$ given by 
 $$g\bullet a  = gag^{-1} + \partial \ell n_G(g).$$
 If $\partial \ell n y = a$ and $g \in G(K)$, then, $\tilde y = g y$ generates the same extension of $K$ as $y$, and satisfies $\partial \ell n (\tilde y) =  g\bullet a.$ 
So, the extension $K(y)/K$, and its Galois group, depend  only on the {\it orbit} $G(K)\bullet a$ of $a$.

\medskip
\noindent
With these points in mind, the main theorem of Kolchin's theory  (see  \cite{vdP-S}, \S 1.4 for the Picard-Vessiot case) can be summarized as follows:

\begin{Theorem} :  i) $Im(\rho) = J_a(C)$, where $J_a/C$   is an algebraic subgroup of $G/C$; 

ii) there is a ``Galois correspondence" : for instance, $K(y)^{J_a(C)} = K$; 

iii) $tr. deg. (K(y)/K) = dim (J_a).$  
\end{Theorem}

The Galois correspondence shows that $J_a$ is connected if and only if $K$ is integrally closed in $K(y)$. To fix the ideas, we shall now  assume that  the base differential field $K$ is algebraically closed, so that all Galois groups over $K$ become connected. We can now describe the Galois goup in a style which will become the  leit-motiv of these notes. See \cite{vdP-S}, I.31.(1)  and I.31.(2) for the Picard-Vessiot case, with a weaker assumption on $K$.

\begin{Theorem} :
  {\it assume that $K$ is algebraically closed. Then, $J_a$ is a minimal algebraic subgroup $J/C$ of $G/C$ such that  $LJ(K)$  meets the orbit of $~a~$ under $G(K)$.}
  
  \end{Theorem}

\noindent
{\it Proof} : consider the set ${\cal H}$ formed by all the algebraic subgroups $H/C$ of $G/C$ such that $G(K)\bullet a \cap LH \neq \emptyset$. Note that this set is stable under conjugation by $G(C)$, since $c(g\bullet a)c^{-1} = (cg)\bullet a $ for any $c \in G(C)$. We will show that any such $H$ contains 
a $G(C)$-conjugate of  $J_a$ - this actually holds for any $K$ - , and that $J_a$ itself lies in this set.

i) Let $H \in {\cal H}$ and let $g \in G(K)$ such that $g\bullet a := \tilde a \in LH(K)$. Since the restriction of $\partial \ell n_G$ to $H$ is $\partial \ell n_H$, which is surjective on $\hat K$-points, there is a solution $ z \in H(\hat K)$ of the equation $\partial \ell n_G (z) = g \bullet a$.  So, $z = \tilde y c = g y c$ with $c \in G(C)$, and the representation $c^{-1} J_ac$ of $\Gamma_a$ attached to $yc$ lies in $H(C)$. 

ii) Consider the $K(y)$-subvariety  $yJ_a$ of $G$. Its set $y J_a(K(y))$ of $K(y)$-points  is Zariski-dense and stable under $\Gamma_a$, and is therefore the set of $K(y)$-points of a $K$-torsor $Z$ under the algebraic group $J_a \otimes_C K$. Since $K$ is algebraically closed, this torsor is trivial, and there exists $u \in Z(K)$ such that $Z = u J_{a/K}$. In particular, $y = u \gamma$ for some $\gamma \in J_a(K(y))$, and $a = \partial \ell n_G (u\gamma) = u \bullet \eta$ where $\eta = \partial \ell n_G \gamma \in LJ_a(K(y))$. Therefore, $u^{-1}Ê\bullet a  = \eta \in LJ_a$ (which must then be $K$-rational), and the  $G(K)$-orbit of $a$ does meet  $LJ_a(K)$.

\bigskip
The above proof  shows that up to $G(C)$-conjucacy,  $\cal H$ contains a {\it unique} minimal element  (we give a more direct proof of this fact in the commutative case in \S 3 below). It does not truly provide an algorithm to compute $J_a$, but it certainly gives upper bounds, which may suffice if one knows  in advance enough elements of  $\Gamma_a$. And of course, a  situation where we can conclude (still with $K$  algebraically closed) is given  by

\medskip

\noindent
{\bf Corollary} :
{\it Assume that $a \in LG(K)$ is {\it non degenerate} : {\it for any proper algebraic subgroup $H \subset G$, the $G(K)$-orbit of $a$ does not meet $LH$}. Then $tr.deg. K(y)/K = dim G$.}

\medskip
When $G$ is commutative, the  action of $G$ on $LG$ reads as : $g \bullet a = a + \partial \ell n_G g$, the set $\partial \ell n_G(G(K))$ is a subgroup of $LG$,  and the theorem  determines the   Galois group   as the smallest algebraic subgroup $H$ of $G$ such that $~a~$ lies in $ \partial \ell n_G(G(K)) + LH$. We derive :

\bigskip
\noindent
1) - {\it  Kolchin's theorem on ${\bf G}_m^n$}: let $y = (y_1, ..., y_n) \in \hat K^{*n}$ such that $\partial y_i /y_i = a_i \in K$ for all $i = 1, ..., n$. Assume that $ \forall  (m_1, .., m_n)  \in \Z^n   \setminus 0, m_1 a_1 + ... + m_n a_n \notin \partial \ell n ( K^*)$. Then,  $tr.deg. K(y)/K = n$. 

Now, if $tr.deg(K/C) = 1$ and if each  $a_i$ is itself of the type $\partial x_i$ for some $x_i \in K$, the condition on the $a_i$'s simply becomes :   $\forall  (m_1, .., m_n)  \in \Z^n   \setminus 0, m_1 x_1 + ... + m_n x_n \notin C$. Indeed, 
$$\partial \ell n(K^*) \cap \partial(K) = \{0\} ~, ~{\rm ~and}~ K\cap \partial^{-1}(0)   = C. \qquad (1)$$  
In the setting of \S 1.2,

\medskip
\centerline{\it we have therefore proved Case (E) of Ax's theorem in the case of tori.}
 
 \medskip
 \noindent
2) - {\it  Ostrowski theorem on $\G_a^n$} : let $x = (x_1, ..., x_n) \in \hat K^n$ such that $\partial x_i = b_i \in K$ for all $i = 1, ..., n$. Assume that $ \forall  (\mu_1, .., \mu_n)  \in C^n \setminus 0,  \mu_1b_1 + ... + \mu_n b_n \notin \partial(K)$. Then,  $tr.deg. K(x)/K = n$. 

Now, if  $tr.deg. (K/C) = 1$ and if each  $b_i$ is itself of the type $\partial \ell n y_i$ for some $y_i \in K^*$, the condition on the $b_i$'s simply becomes :   $\forall  (m_1, .., m_n)  \in \Z^n   \setminus 0, y_1^{m_1} ... y_n^{m_n}  \notin C^*$, i.e. $m_1 x_1 + ... + m_n x_n \notin C$. Indeed, the natural map 
$$\iota :  C \otimes_\Z \partial \ell n(K^*) \rightarrow K ~{\rm is ~injective}~,  Im(\iota) \cap \partial(K) = \{0\}  ~, ~{\rm  and}~ K^* \cap \partial \ell n^{-1}(0) =  C^*. \qquad (2)$$
In the setting of \S 1.2,

\medskip
\centerline{\it we have therefore proved Case (L) of Ax's theorem in the case of tori.}

\bigskip
We refer to \cite{BM}, \S 6,  for a similar treatment  of Cases (E) and (L) of Ax's theorem on general commutative algebraic groups $G$ defined over $C$. The required analogues $(1^*), (2^*)$ of the displayed Formulae (1),(2) are discussed in \S 4.1 below.

\section{$D$-groups and Pillay's theory}

\subsection{General setting}
 In a series of papers started in 1997 \cite{Pi1}, A. Pillay extended Kolchin's theory to the context of algebraic $D$-groups. These groups, which, after the work of P. Cassidy (see \cite{B-C}, \S 2),  had been considered by Buium \cite{Bu} with  an eye towards  the  functional analogue of the Manin-Mumford and Mordell-Lang conjectures, are defined over a differential field $(K, \partial)$, and usually not over its field of constants $C$ (but even then, new phenomena can occur, cf. \S 3.2). For the sake of simplicity, we will restrict the presentation of Pillay's theory to the case of  {\it commutative} $D$-groups, with an additive notation. Again, we suppose that $K^\partial = C$ is algebraically closed of characteristic $0$, and we fix a differential closure $\hat K$ of $K$.  

So,  let 
$G/K$ be a connected commutative algebraic 
group over $K$.   We 
 {\bf assume} that $\partial$ extends to a derivation
 $D_\partial : {\cal O}_{G_K}Ê\rightarrow{\cal O}_{G_K}$  of  the structure sheaf ${\cal O}_{G_K}$   compatible with its structure of Hopf algebra. (The derivation $D^0_\partial$ considered in \S 2 when $G$ is defined over $C$ did satisfy this property.) We {\bf fix} such an extension $D_\partial$, and say that  $G$, tacitly equipped with $D_\partial$, {\it  is a $D$-group}. We can then define the logarithmic derivative of $G$  in exactly the same way as before : for any $y \in G(K)$, we have the two derivations $D_{\partial, y}$, $\delta y$ of the local ring ${\cal O}_{G,y}$, and $\partial \ell n_G(y) \in LG(K)$ is their difference $\delta y - D_{\partial, y}$, pulled back from $T_yG$ to $LG$ via the canonical splitting of the tangent bundle $TG \simeq G \times LG$. Again, this extends over any differential extension $\cal K$ of $K$, and we set $G^\partial({\cal K}) = \{ y \in G({\cal K}), \partial \ell n_G (y) = 0\}$. By \cite{Pi}, 2.5, $\partial \ell n_G$ is surjective at the level of $\hat K$-points.

  \medskip

Given $a \in LG(K)$, Pillay's theory studies the differential extension $K(y)/K$, where $y$ is a  solution in $G(\hat K)$ of 
$$\partial \ell n_G(y) = a,$$
and its Galois group $\Gamma_a = Aut_\partial(K(y)/K) $. To avoid
\medskip
 ``new constants", we now request that 
 \medskip
\centerline{Ê$G$ is ``$K$-large" : $G^\partial(\hat K) = G^\partial(K)$;}
\medskip
\noindent
Notice that this hypothesis was automatically  satisfied in the case of \S 2. 

\medskip

The fact that $D_\partial$ respects the group structure of $G$ is again reflected by a cocycle identity which, in our commutative situation, becomes :
$$\forall u, v \in G , \partial \ell n_G(uv) = \partial \ell n_G u + \partial \ell n_G v.$$
Consequently,

i) two solutions differ by an element of $G^\partial(\hat K) = G^\partial(K)$, hence  define the same extension of $K$, and the Galois group $\Gamma_a$ comes with a canonical embedding $\xi$ into $G^\partial(K)$ : 
$$\forall \sigma \in \Gamma_a~, ~\sigma y - y  :=  \xi(\sigma)  \in G^\partial(\hat K) =G^\partial(K);$$

ii) the extension $K(y)/K$ and its Galois group depend only on the image of $~a~$ in the group
$$Coker(\partial \ell n_G, K) := LG(K) /\partial \ell n_G (G(K)).$$
We extract from the main theorem of Pillay's theory (see \cite{Pi}, \S 3) :

\begin{Theorem}:  assume that the $D$-group $G$ is $K$-large. Then :

{i) \it  $Im(\xi) = N_a^\partial(K)$, where $N_a/K$   is a   $D$-subgroup of $G/K$;
  
  ii)  there is a ``Galois correspondence" : for instance, $K(y)^{N^\partial_a(K)} = K$; 
  
  iii)  $tr. deg. (K(y)/K) = dim (N_a).$}
  \end{Theorem}

 Here, a $D$-subgroup $H$ of $G$ is an algebraic subgroup of $G$ whose ideal sheaf ${\cal I}_H \subset {\cal O}_G$ is stable under  $D_\partial$, and the $D$-structure of $H$ is given by the derivation induced by $D_\partial$ on ${\cal O}_H$. In particular, $(\partial \ell n_G)_{ | H} = \partial \ell n_H$, $H$ is  $K$-large, and $\overline G = G/H$ acquires a natural structure of $D$-group, which is $K$-large as well. (Notice that the sequence $0 \rightarrow H^\partial(\hat K) \rightarrow G^\partial (\hat K) \rightarrow \overline G^\partial (\hat K) \rightarrow 0$ is exact, since $\partial \ell n_H$ is surjective on $\hat K$-points.) Still assuming $K$-largeness, we have :

\begin{Theorem} :  the identity component of   $N_a$ is the smallest $D$-subgroup $H/K$ of $G/K$ such that 
$a$ lies in $LH + \Q. \partial \ell n_G G(K).$
\end{Theorem}
 \noindent
 Suppose for a moment that $K$ is algebraically closed. Then  $N_a$ is connected and $ \partial \ell n_G G(K)$ is already a $\Q$-vector space, so the theorem acquires the same   form as Theorem 2.2. Here, the commutativity assumption allows to drop any requirement on  $K$.

\medskip
\noindent
{\it Proof} :  as in \S 2, consider the set ${\cal H}$ formed by all the $D$-groups $H/K$ of $G/K$ such that $ a \in LH(K) + \Q.\partial \ell n_GG(K)$. Then,

i)  for any $H$ in $\cal H$,  the equation corresponding to some multiple of $a$ has a solution in $H(\hat K)$, so the connected component of $Im(\xi)$ lies in $H$.

ii) the image $\overline y$ of $y$ in $\overline G = G/N_a$ is stable under $\Gamma_a$, hence $K$ rational. Analysing commutative algebraic groups, we see that  $G(K)$ projects onto a subgroup of finite index of $\overline G(K)$, so $m\;y = u +  \gamma$ for some $m > 0, u \in G(K), \gamma \in N_a(K(y))$,  and $m\;a = \partial \ell n_G(u) + \eta$, where $\eta = \partial \ell n_{N_a}(\gamma)$ is a (necessarily $K$-rational) point of $LN_a$. 

\medskip
As promised in \S 2 , we now give a direct proof  (independent of Galois theory) that $\cal H$ admits a unique minimal element.  This fact depends crucially on the hypothesis of $K$-largeness of $G$, which implies that for any $D$-subgroup $H$ of $G$,  $ LH(K) \cap \partial \ell n_G(G(K)) = \partial \ell n_H(H(K))$, i.e.  that the natural map 
$Coker(\partial \ell n_H , K) \rightarrow Coker(\partial \ell n_G, K)$ is injective. The claim   then easily follows.  When $G$ is not $K$-large, the snake lemma merely says that the sequence 
$$G^\partial(K)/H^\partial(K) \hookrightarrow (G/H)^\partial(K) \rightarrow Coker(\partial \ell n_H , K) \rightarrow Coker(\partial \ell n_G, K)$$
is exact.

 \medskip
 Anyway, the outcome in the $K$-large case is exactly the same as in the constant one. For instance,  assume that $a$ is {\it non degenerate} : {\it for any proper algebraic $D$-subgroup $H \subset G$,  $a + \Q.\partial \ell n_G G(K)$ does not meet $LH(K)$}. Then $N_a = G$, i.e.  Kolchin-Ostrowski still holds true. In \S 4, we will find  conditions on $a$ ensuring its non-degeneracy, and thereby obtain  ``non-constant" analogues of the $(E)$ and $ (L)$ statements. But more work is required before we get there, because  the $K$-largeness hypothesis is seldom satisfied in these non-constant situations. A good example (and a way to overcome the difficulty) is given by the following case.

\subsection{The case of $D$-modules}

Let $B \in \frak{gl}_n(K)$. In \S 2, we used Kolchin's view-point to describe the Picard-Vessiot extension attached to the homogeneous equation $\partial y = B y , y \in \hat K^n$.  Given $a \in K^n$, we now study  the inhomogeneous equation 
$$\partial y = B y + a.$$
Setting $A =  \left( \begin{array}{ccc}
 B& a \\
0 &0
 \end{array} \right) \in \frak{gl}_{n+1}(K)$,
 and $Y =  \left( \begin{array}{ccc}
y \\
1
 \end{array} \right) \in \hat K^{n+1}$, this equation is essentially the same as $\partial Y = A Y$, and can indeed be treated by pure Picard-Vessiot means, as was done in \cite{BS}, \cite{B1}. But we will now describe it from the point of view of Pillay's theory, starting with a vectorial group $G := V \simeq K^n$ over $K$.
 
 \medskip
Given a basis of $V$ over $K$, the choice of derivation $D_\partial$ on the affine ring $K[Sym(V^*)] \simeq  K[X_1, ..., X_n]$  of $V$ inducing $\partial$ on $K$ and respecting the group structure of $V$ is tantamount to the choice of a matrix $B \in Mat_{nn}(K)$ such that $D_\partial(X_1, ..., X_n) = (X_1, ..., X_n)B$, i.e. of a $D$-module structure on $V$. The associated logarithmic derivative, denoted in this vectorial case\footnote{~In this paragraph, we  keep to the notation $\partial \ell n_V$ to ease the comparison with the previous cases, but we will state the final result with $\nabla_V$ to prepare for the study of the $D$-module $L\tilde A$ in \S 4.2.} by
$$\nabla_V = \partial \ell n_V : V \rightarrow LV \simeq V~,$$
 is then given by 
$$V(K)  \simeq K^n \ni y \mapsto \partial \ell n_V(y) = \partial y - B y \in LV(K) \simeq K^n,$$
where $\partial ~^t(y_1, ..., y_n) = ~^t(\partial y_1, ...., \partial y_n)$.  Then,
$$V^\partial(\hat K) = \{y \in \hat K^n, \partial y = By\}$$
(which we will abusively denote by $V^\partial$) is the $C$-space of solutions of a linear equation, so that $V$ is usually not $K$-large. In order to apply Pillay's theory, we extend $K$ to the Picard-Vessiot extension $\K := K_V := K(V^\partial(\hat K))$. By definition, $V_\K = V \otimes_K \K$ is now clearly $\K$-large.

\medskip
Given a $K$-rational point $a \in LV \simeq V$, we consider the equation
$$\partial \ell n_V y = a ~, ~i.e. ~\partial y - By = a.$$
By Theorem 3.2 and the $\Q$-divisibility of rational points in vectorial groups, its Galois group $\Gamma_a = Aut_\partial({\K}(y)/{\K}) $ is a $C$-subspace of $V^\partial$ of the form  $N_a^\partial =  N_a^\partial(\hat K) = (N_a)^\partial({\K})$, where $N_a \simeq LN_a$ is the smallest $\K$-vector subspace of $V_\K$  stable under $\partial \ell n_V$,  such that $a \in \partial \ell n_V(V({\K})) + N_a$. Let us try to turn this into a more manageable description.

\medskip 
 A first remark is that $N_a$ is actually defined over $K$. Indeed, let ${\cal H}_\K$ be the set of all $\K$-subspaces of $V_\K$ satisfying the property above. Since $a$ is defined over $K$, the set ${\cal H}_\K$ is stable under  the action of $J(C) := Aut_\partial(\K/K)$, and  its unique minimal element $N_a$ too is stable under  $J(C)$.  By Picard-Vessiot theory, $N_a$ must then be defined over $K$. Another  proof of this fact will be given presently, cf. Proof (i) below.
 
In these conditions, it is tempting to consider the set ${\cal H}_K$ of all $K$-vector subspaces $N$  of $V$  stable under $\partial \ell n_V$,  such that $a \in \partial \ell n_V(V(K)) + N(K)$. Since $  {\cal H}_K \subset {\cal H}_\K $, the Galois group $N_a$ is contained in all its elements. But in general, $N_a$ {\bf will not belong} to ${\cal H}_K$, and a priori, one cannot even speak of the smallest element of ${\cal H}_K$. See \cite{B1} for some counterexamples.

As noticed in \cite{BS}, \cite{B1}, there is however a case where $N_a$ does belong to ${\cal H}_K$, and therefore becomes its smallest element, viz. when the differential system $\partial y = By$ can be split over $K$ into a direct sum of irreducible systems, i.e.  when the $D$-module structure on $V$ attached to $D_\partial$ is semi-simple over $K$. Then, any $D$-submodule  $H/K$ of $V$ admits a $D$-submodule complement  over $K$,  and this yields the injectivity of $Coker(\partial \ell n_H , K) \rightarrow Coker(\partial \ell n_V, K)$. So we can already say that ${\cal H}_K$  has a unique minimal element. Furthermore, as shown  in Proof (ii) below, given any  quotient $\overline V = V/H$ of    such a semi-simple $V$, the natural map 
$$Coker(\partial \ell n_{\overline V}, K) \rightarrow Coker(\partial \ell n_{\overline V}, \K)$$
is injective. So, for $H$ in ${\cal H}_\K$,  the  image $\overline a$ of $a$ in $L\overline V (K) $   lies in $\partial \ell n_{\overline V}(\overline V(\K))$ only if it already lies in $\partial \ell n_{\overline V}( \overline V(K))$, in which case $a + \partial \ell n_V (V(K))$ meets $LH$. Consequently, any $H$ in ${\cal H}_{\K}$ contains an element of ${\cal H}_K$, the minimal elements of ${\cal H}_{\K}$  and of  ${\cal H}_K$ coincide, and we derive as in 
 \cite{BS}, \cite{B1} (see also \cite{Ch}, Lemme 2.2.10 and Thm. 2.2.5) :

\begin{Theorem} : assume that the $D$-module $(V/K, \nabla_V)$ is semi-simple. Then, the Galois group $Aut_\partial({\K}(y)/{\K})$ is  $N_a^\partial$, where $N_a$  is the smallest $D$-submodule $N$ of $~V$ defined over $K$ such that $~a~$ lies in $  \nabla_V (V(K)) + N$.
\end{Theorem}

\noindent
{\it Descent proofs} (as promised above) : (i) Firstly, the other   proof that under  no hypothesis on $V$, $N_a$ is always defined over $K$. Consider the tower of extensions, all Picard-Vessiot {\it over $K$} in view of the system $\partial Y = A Y$ mentioned at the beginning :

\begin{equation*}
\begin{array}{cccccc}
&&\K(y) & & \xi  \\
&& \mid &  \} N^\partial_a  & \hookrightarrow & V^\partial\\
&\Gamma \{&\K  & &  \rho \\
&&\mid & \} J & \hookrightarrow & GL(V^{\partial})  \\
&&K  \\
 \end{array}
 \end{equation*}
The full Galois group $Aut_\partial(\K(y)/K)$ is of the form $\Gamma(C)$ for some algebraic subgroup $\Gamma/C$ of $GL_{n+1}$. Since  $\K/K$ is P-V, $N^\partial_a$ is a normal subroup of $\Gamma $, with quotient $J$ naturally acting on $V^\partial$ by a $C$-rational representation $\rho : J \rightarrow GL(V^\partial) \simeq GL_{n/C}$ of the type described in \S 2. Since $N_a$ is abelian, $J$ also acts on $N_a^\partial$ by conjugation, and  a familiar  computation from affine geometry shows that the homomorphism $\xi$ commutes with these actions of $J$ :
$$\forall \sigma \in N_a^\partial, \tau \in J~, ~ \xi(\tau \sigma \tau^{-1}) = \rho ( \tau) \big(  \xi(\sigma)\big).$$
Therefore, $\xi$ identifies $N^\partial_a$ with a $J$-submodule of $V^\partial$, and by the standard P-V theory of \S 2, $N_a$ must then be a $D$-submodule of $V$ defined over $K$.

\medskip
The $J$-equivariance property on which this proof  is based can be viewed as a first instance of the  arguments of Galois descent  alluded to in the introduction of the paper. Indeed, $\xi$ is the restriction to $N_a^\partial$ of a $C$-rational cocycle $\hat \xi_a   : \Gamma \rightarrow V^\partial$, defined by the same formula $\hat \xi_a(\tau) = \tau y - y$, whose class in $H^1(\Gamma, V^\partial)$ depends only on the image $\tilde a$ of $~a~$ in  $Coker(\partial \ell n_V, K)$. More precisely, the map 
$$\Xi_K :  Coker(\partial \ell n_V, K) \rightarrow H^1 (\Gamma, V^\partial) : \tilde a \mapsto \Xi_K(\tilde a) = {\rm{~class ~of ~}} \hat{\xi_a}$$
is a   group embedding. Now, it is a well-known feature of group cohomology that  $\hat \xi_a$, restricted to any normal subgroup $\rm N$ of $\Gamma$, provides a $\Gamma/{\rm N}$-invariant cohomology class in $H^1({\rm N}, V^\partial)$, cf. \cite{Se}, I.2.6.b.

  \medskip
  (ii) We now assume that $V$ is semi-simple. Any quotient $\overline V$ is then also semi-simple, and we must prove that the kernel of the map  $Coker(\partial \ell n_{\overline V}, K) \rightarrow Coker(\partial \ell n_{\overline V}, \K)$ vanishes. But via the map $\Xi_K$ attached to $\overline V$ and the analogous map $\Xi_\K$ at the level of $\K$, this kernel injects into the kernel of the restriction map
  $$H^1(\Gamma, \overline V^\partial) \rightarrow H^1(N_a^\partial, \overline V^\partial)^J = Hom_J(N_a^\partial, \overline V^\partial).$$
The latter kernel identifies with $H^1(J, \overline V^\partial)$ by the inflation map \cite{Se}, {\it loc. cit.}. Finally, the  faithful representation $V^\partial$ of $J$ is  completely reducible,  so $J$  is a reductive group and  $H^1(J, \overline V^\partial) = 0$. (Notice that  all cocycles appearing here are  continuous for the Zariski topology. Even when $J = PSL_2$ and $C = \C$, I do  not know if the abstract cohomology  groups of  the abstract group $J(C)$ would also vanish.)

\bigskip
In the next section, we show that a similar descent argument applies to various $D$-groups attached to abelian varieties. Obstructions occur in the case of   semi-abelian varieties, and we refer to  \cite{An}, \cite{B-P}, \cite{BM} for examples and counterexamples they lead to.

\section{Abelian varieties}

From now on, we fix a smooth  irreducible curve $S$ over the field $\C$ of complex numbers,  and a non-zero vector field $\partial \in H^0(S, TS)$ on $S$, which we identify with a  derivation on the function field $K = \C(S)$, with $K^\partial = \C$. We denote by $\overline K \subset \hat K$ its algebraic closure. In the sequel, we may tacitly restrict $S$ to a non empty open subset, or more generally to a finite cover, and still denote by $S$ the resulting curve. We consider an abelian variety $A/K$, extended to an abelian scheme $\  {\cal A}/ S$, with relative  Lie algebra $L{\cal A}/S$.  We assume that   the largest abelian variety $A_0/\C$, isomorphic over $\overline K$ to an abelian subvariety of $A$, is embedded in $A$, and call it the  constant part, or $\C$-trace, of $A$.  

\medskip

The exponential maps of the various fibers of ${\cal A}/S$ provide a morphism $exp_A : L{\cal A}^{an}  \rightarrow {\cal A}^{an}$ of analytic sheaves over the Riemann surface $S^{an}$. For a local section $\overline x$ of $L{\cal A}^{an} $, we denote by $\overline y = exp_A(\overline x)$ its image in ${\cal A}^{an}$. We say that 

\medskip
   {\it $\overline x \in LA$ is non-degenerate}  if for any proper abelian subvariety $B$ of $A$, $\overline x$ does not lie in $LB + LA_0(\C)$;

{\it $\overline y \in A$ is non-degenerate} if for any proper abelian subvariety $B$ of $A$, no non-zero multiple of $\overline y$  lies in $B +  A_0(\C)$. Then:

\begin{Theorem} :  let $A/K$ be an abelian variety, with $\C$-trace $A_0$, and let $\overline x \in LA$, with image $ \overline y = exp_A(\overline x) \in A$. Assume that
 
 (L)  $\overline y \in A$ is $K$-rational and non-degenerate; then,  $tr.deg. K(\overline x)/K = dim   A$;
 
 (E) $\overline x \in LA$ is  $K$-rational and non-degenerate; then,  $tr.deg. K(\overline y)/K = dim   A$.
 \end{Theorem}
 \noindent
In case (L), one cannot replace the non-degeneracy hypothesis on $\overline y$ by the weaker one on $\overline x$,   because the  periods of $exp_A$ may satisfy non linear algebraic relations{\footnote{~This problem can be addressed through the study of the group $J$ appearing below. Applications of the type discussed in \cite{MZ} may require such sharpenings.}. Similarly, under the mere hypothesis that $\overline x$ is non-degenerate, the analogue of the Ax-Schanuel theorem  (that $tr. deg.(K(\overline x, \overline y)/K) \geq dim A$) would not hold, but  one can conjecture that it does as soon as $\overline y$ is non-degenerate.

The theorem reflects  the existence of large Galois groups attached to $A/K, \overline x, \overline y$. But before we can speak of Galois groups, we need a $D$-group.

\subsection{From abelian varieties to $D$-groups}

In general, given a commutative algebraic goup $G/K$, the set of all possible extensions of $\partial$ to a derivation $D_\partial$   on ${\cal O}_G$ is empty, or is an affine space under the space of sections of the tangent bundle $TG$. It therefore corresponds to a class  $\kappa(G/K, \partial)$ in $H^1(G, TG)$ and is non-empty  only if this class vanishes. When $G = A$ is proper, $\kappa(A/K, \partial) $ is the image of $\partial$ under the Kodaira--Spencer map attached to ${\cal A}/S$ at the generic point of $S$, and is known to vanish if and only if  ${\cal A}/S$ is isoconstant, i.e.    $A \simeq A_0$ is isomorphic over $\overline K$ to an abelian variety  defined over $\C$. So, {\it a non-isoconstant abelian variety $A/K$ admits no $D$-group structure}.

\medskip
To overcome this difficulty (see   \cite{Bu}), we introduce the universal extension $\tilde A/K$ of $A$. This is, in the category of algebraic goups, an extension
$$ 0 \rightarrow W_A \rightarrow \tilde A \rightarrow A \rightarrow 0$$
of $A$ by a vectorial group $W_A/K$, dual to $H^1(A, {\cal O}_A)$. In particular, $\tilde A$ has dimension $2 dimA$, and  in fact, its Lie algebra 
$L\tilde A$ is dual to the de Rham cohomology group $H^1_{dR}(A/K)$ of $A/K$. Now, the latter admits a natural connection (Gauss-Manin), whose  dual $\nabla_{L\tilde A}$, contracted with $\partial$, provides a {\it $D$-module structure on $L\tilde A$}. Finally,  $\nabla_{L\tilde A}$ can be ``integrated" into a {\it $D$-group structure on $\tilde A$}, which is actually unique. We point out for later use that in the usual identification of a vectorial group with its Lie algebra, {\it the restrictions of $\nabla_{L\tilde A}$ and $\partial \ell n_{\tilde A}$ to $W_A \simeq LW_A$ coincide}, so that statements comparing   their values often factor through the quotients $A, LA$. But unless $A/K$ is isoconstant, $W_A$ is not a $D$-subgroup  or a $D$-submodule of $\tilde A, L\tilde A$. Another property to keep in mind in what follows is that  any algebraic subgroup of $\tilde A$ projecting onto $A$ must fill up $\tilde A$;  by Poincar\'e and the fonctoriality of universal extensions, this implies that any connected $D$-subgroup  of $\tilde A$ is of the form $\tilde B + W$, where $B$ is an abelian subvariety of $A$ and $W$ is a $D$-submodule of $L\tilde A$ contained in $W_A$. 

\medskip
We now give an analytic description of the  logarithmic derivative $\partial \ell n_{\tilde A}$ of the $D$-group structure of $\tilde A$,  in the style of \S 1.2.  Extend $\tilde A/K$ to a group scheme $\tilde {\cal A}/S$, and consider the exact sequence of analytic sheaves over  $S^{an}$  given by the exponential morphism :
$$0 \rightarrow \Omega_{\tilde {\cal A}} \rightarrow L\tilde {\cal A}^{an} \rightarrow^{exp_{\tilde A}}  \tilde {\cal A}^{an} \rightarrow 0.$$
Its kernel $\Omega_{\tilde {\cal A}}$ is the $\Z_{S^{an}}$- local system of  periods of $\tilde {\cal A}$, and by the analytic definition of the Gauss-Manin connection,  these generate over $\C_{S^{an}}$ the  space $(L\tilde A)^\partial$ of horizontal sections of the connexion  $\nabla_{L\tilde A}$. Therefore, given $y  \in \tilde A(K)$, extended to a section $y(z)$ of $\tilde {\cal A}/S$, and any local choice $\ell n_{\tilde A} {y}(z)$ of an inverse of $ y$ under $exp_{\tilde A}$,  $\nabla_{L\tilde A} \circ \ell n_{\tilde A} {y(z)}$   extends to a well-defined section of $L\tilde A$ over $S^{an}$, actually with moderate growth at infinity, hence over $S$. Finally, as shown in \cite{B-P}, Appendix H, the resulting point in $L\tilde A(K)$ 
coincides with the logarithmic derivative  $\partial \ell n_{\tilde A} (y)$ of $y$. In brief, $\partial \ell n_{\tilde A} = \nabla_{L\tilde A}Ê\circ \ell n_{\tilde A}$, and the analytic relation $y = exp_{\tilde A}(x)$ implies
$$ \partial \ell n_{\tilde A} y = \nabla_{L\tilde A} x.$$
Conversely, when $x$ and $y$ have an analytic meaning, this differential relation is equivalent to the existence of 	a local horizontal section $\xi \in   L\tilde A, \nabla_{L\tilde A}(\xi) = 0$, of $ \nabla_{L\tilde A}$ such that $y = exp_{\tilde A}(x-\xi)$. But contrary to the  situation described in \S 1.2, the point $\xi$ is not necessarily in the constant part $L\tilde A_0(\C)$; neither is $\eta = exp_{\tilde A}(\xi)$, which, in the elliptic case, can be  computed in terms of   solutions of Picard  type of certain Painlev\'e VI equations. 

\begin{Theorem}Ê:  let $A/K$ be an abelian variety, with $\C$-trace $A_0$, and let $x  \in  L\tilde A (\hat K), y \in \tilde A(\hat K)$ satisfy the differential relation $ \partial \ell n_{\tilde A} y = \nabla_{L\tilde A} x$. Assume that  
  
{\rm  (L) $=$  \cite{An}, Theorem 3} :    $y  \in  \tilde A(K)$, and its image $~\overline y$ on   $ A(K)$ is non-degenerate; then  $tr.deg. K(x)/K = 2\; dim  A$; 
 
{\rm  (E) $=$  \cite {B-P}, Theorem 1.4 }  :      $x \in L\tilde A(K)$, and its image $\overline x \in LA(K)$ is non-degenerate; then, $tr.deg. K(y)/K = 2 \;dim  A$. 
\end{Theorem}
 
 \medskip
 The theorem announced at the beginning of the section easily follow, by the functoriality of the exponential morphisms,   the fact that any $K$-rational point of $LA, A$ can be lifted to a $K$-rational point of $L\tilde A, \tilde A$, and a dimension count. And the present theorem  will clearly follow from Theorems 4.3 and 4.4 below, since its non-degeneracy hypotheses mean, in their notations, that $B = A$, and transcendence degrees, controled by the dimension of Galois groups in view of Theorem  3.1.iii, cannot decrease when we go from $\K$ to $K$. Finally, we can  replace $K$ by its algebraic closure  in these statements, so, we henceforth assume that $K = \overline K$ is an {\it  algebraically closed field} of transcendence degree 1 over $\C$.

 \medskip

We  now collect  some facts about the $D$-groups $\tilde A, L\tilde A$, in particular with respect to their (non)-$K$-largeness, referring to \cite{B-P} and  \cite{BM} for their proofs. Suffices to say here that those concerning $L\tilde A$ are based on deep theorems on variations of Hodge structures, cf. \cite{De}, while those concerning $\tilde A$ depend on a combination of the latter with similarly deep theorems of model theory, cf. \cite{M-P}, and with Manin's theorem below.  

\medskip
\noindent
{\it Fact 1} : $(L\tilde A)^\partial(K) = L\tilde A_0(\C)$. So, $L\tilde A$ is $K$-large only if $A/K$ is isoconstant. In general, we denote by $\K = K((L\tilde A)^\partial)$ the extension  generated by the solutions of $\nabla_{L\tilde A} (\xi) = 0$. This is a Picard-Vessiot extension of $K$, and we call $J$ its  Galois group.

\noindent
{\it Fact 2} : the $D$-group $\tilde A/K$  is often $K$-large (for instance, as soon as the Kodaira-Spencer map $\kappa(A/K, \partial)$ has maximal rank), but not always. However, the field $K({\tilde A}^\partial)$ generated over $K$ by the elements of $\tilde A^\partial(\hat K)$ is always contained in the above Picard-Vessiot extension $\K$ of $K$. So, viewed over $\K$, both  $D$-groups $L\tilde A$ (by definition) and $\tilde A$  are $\K$-large.

\noindent
{\it Fact 3} : $L\tilde A$ is a semi-simple $D$-module over $K$. So,  the (connected) algebraic group  $J$ is reductive, and actually, $J$ is even  semi-simple  ($J = [J,J]$). 

\noindent
{\it Fact 4} : let $U_A$ be the maximal $D$-subgroup of $L\tilde A$ contained in $W_A$. Then,  the  quotient $\tilde A/U_A$ is $K$-large.

\medskip

We will also need the analogues of the displayed formula  (1) and its sharper version (2) given at the end of \S 2 in the case of tori.   These are provided by sharper and sharper versions of Manin's kernel theorem, whose simplest case reads as follows : $\partial \ell n_{\tilde A}(\tilde A(K) ) \cap \nabla_{L\tilde A}(L\tilde A(K)) = \nabla_{L\tilde A}(W_A(K))$, and more precisely :

\medskip

\centerline{$y \in \tilde A(K) , \partial \ell n_{\tilde A}(y) \in \nabla_{L\tilde A}(L\tilde A(K)) \Rightarrow$ a multiple of $\overline y$ lies in $   A_0(\C). \qquad (1')$}

\medskip
\noindent
(By ``a multiple", we mean a multiple by a {\it non-zero} integer; this convention will be used throughout the rest of the text.)   This was extended by Chai to quotients of $\tilde A$ by $D$-subgroups of $\tilde A$, as follows : let $H$ be a  $D$-subgroup of $\tilde A$, with image $\overline H$ in $A$; then

\medskip
\centerline{$y \in \tilde A(K) , \partial \ell n_{\tilde A}(y) \in \nabla_{L\tilde A}(L\tilde A(K)) + LH  \Rightarrow$ a multiple of $\overline y$ lies in  $ A_0(\C) + \overline H, \quad (1^*)$}

\medskip
\noindent 
and further extended to quotients of $L\tilde A$ by $D$-submodules, as follows :  assume that $A/K$ is a {\it simple}   abelian variety, and let $N$ be any proper $D$-submodule of $L\tilde A$; then 

\medskip
\centerline{$y \in \tilde A(K) , \partial \ell n_{\tilde A}(y) \in \nabla_{L\tilde A}(L\tilde A(K)) + N \Rightarrow$ a multiple of $\overline y$ lies in  $    A_0(\C) \quad (2').$}

\medskip
See  \cite{B2}, \S 3, for a proof of (2'), based on Andr\'e's normality theorem \cite{An}, Thm. 2, from which yet another sharpening can be deduced, as follows. Recall the definition of non-degeneracy of $\overline y \in A$ given at the beginning of the section,  let $A/K$ be any abelian variety, and let $N$ be any proper $D$-submodule of $L\tilde A$; then 
 
 \smallskip
\centerline {$y \in \tilde A(K) , \partial \ell n_{\tilde A}(y) \in \nabla_{L\tilde A}(L\tilde A(K)) + N \Rightarrow \overline y$   is degenerate in   $A$.  \qquad ($2^*$) }

\medskip
\noindent
Indeed,  the semi-simplicity of $L\tilde A$ allows us to speak of the smallest $D$-submodule $N_0$ satisfying this property. By the normality theorem, there exists a $D$-subgroup $H$ of $\tilde A$ such that $N_0 = LH$, and since $N$ is proper, $\overline H \neq A$. The conclusion now follows from ($1^*$). 

\medskip
The reader will notice that Formulae (1) and (2) of \S2 on a torus $\G_m^n$   are the exact analogues of  the abelian Formulae ($1^*$) and ($2^*$).

 \subsection{Abelian logarithms}
 
 Let $A/K$ be an abelian variety, with $\C$-trace $A_0$,  and let $y \in \tilde A(K)$ with image $\overline y \in A(K)$.  Recall the notation $\K = K((L\tilde A)^\partial)$   from Fact 1, the assumption that $K$ is algebraically closed 
  and the convention that a multiple means a multiple by a non-zero integer.
 
We consider the differential system in the unknown $x \in L\tilde A(\hat K)$ :
  $$\nabla_{L\tilde A} x = a  ~{\rm where~} ~ a =  \partial \ell n_{\tilde A} y \in L\tilde A(K). \qquad (*).$$
 Its Galois group over $\K$  is well-defined, since the $D$-module $L\tilde A$ is $\K$-large, and is the normal subgroup $N_a^\partial$ of $\Gamma$ in the diagram of \S 3.2, which here reads :
 
\begin{equation*}
\begin{array}{cccccc}
&&\K(x) & & \xi  \\
&& \mid &  \} N^\partial_a  & \hookrightarrow & (L\tilde A)^\partial\\
&\Gamma \{&\K  & &  \rho \\
&&\mid & \} J & \hookrightarrow & GL((L\tilde A)^{\partial})  \\
&&K &,  \\
 \end{array}
 \end{equation*}
where $\xi(\sigma) = \sigma(x) -x$, and $J$ acts $\C$-rationally on $(L\tilde A)^\partial$. 
 
\begin{Theorem}Ê: let $y \in \tilde A(K)$, and let $B$ be the smallest abelian subvariety of $A$ such that  some multiple of  $\overline y$   lies in $B+ A_0(\C)$. Then, the Galois group  $Aut_\partial (\K(x)/\K)$ of $(*)$ over $\K$ is equal to $(L\tilde B)^\partial$, and in particular, $tr.deg.(\K(x)/\K) = 2 dim B$.
\end{Theorem}

 \noindent
 {\it Proof} : by Fact 3 and  Theorem  3.3,  this Galois group is $N_a^\partial$, where $N_a$ is the smallest $D$-submodule $N/K$  of $L\tilde A$ such that
$a := \partial \ell n_{\tilde A} (y) \in N + \nabla_{L\tilde A}(L\tilde A (K))$.

\medskip

Fix a lift $y' \in \tilde B(K)$ of  $\overline y$. Then,  $u = y-y' $ lies in $W_A(K)$, where $\nabla_{L\tilde A}$ and $\partial \ell n_{\tilde A}$ coincide. Therefore, $\partial \ell n_{\tilde A}(y)  = -\partial \ell n_{\tilde B}(y') + \nabla_{L\tilde A}(u)   \in L\tilde B + \nabla_{L\tilde A}(L\tilde A(K))$, and $N_a $ is contained in $L\tilde B$.

\medskip
To prove the converse inclusion,  we set $x' = x-u$, and note that by Fact 1, the Galois group of ($*$) over $\K$ is the same as that of $\nabla_{L\tilde B} x' = a'  ~{\rm where~} ~ a' =  \partial \ell n_{\tilde B} y' \in L\tilde B(K)$. Furthermore, the Picard-Vessiot extensions $\K$ and $\K_B(x')$  of $\K_B = K((L\tilde B)^\partial)$ are linearly disjoint, since a normal subgroup of the reductive group $J$ cannot admit a non-zero vectorial quotient. So $N_a$ is in fact  the smallest $D$-submodule $N/K$  of $L\tilde B$ such that
$  \partial \ell n_{\tilde B} (y') \in N + \nabla_{L\tilde B}(L\tilde B (K))$.  Since the image $\overline y'$ of $y'$ in the abelian variety $B$ is by definition non-degenerate, the strong form ($2^*$) of Manin--Chai, applied to $\tilde B$, now  implies that $N_a$ fills up $L\tilde B$.  

 \subsection{Abelian exponentials}
 
 Let $A/K$ be an abelian variety, with $\C$-trace $A_0$.  Recall the notation $\K = K((L\tilde A)^\partial)$,   and let now $x \in L\tilde A(K)$ with image $\overline x \in LA(K)$.  
 We consider the differential system in the unknown $y \in \tilde A(\hat K)$ :
  $$\partial \ell n_{\tilde A} y= a  ~{\rm where~} ~ a =  \nabla_{L\tilde A} x  \in L\tilde A(K). \qquad (**).$$
 Its Galois group over $\K$ is well-defined, since the $D$-group $ \tilde A$ is $\K$-large, and is represented by the top level of the tower of extensions :
 \begin{equation*}
\begin{array}{cccccc}
&&\K(y) & & \xi  \\
&& \mid &  \} H^\partial_a  & \hookrightarrow & \tilde A^\partial \\
&{\rm no }~ \Gamma ~! &\K  & &  \rho \\
&&\mid & \} J(C) & \rightarrow &  Aut(\tilde A^{\partial})  \\
&&K &,  \\
 \end{array}
 \end{equation*}
where $\xi(\sigma) = \sigma(y) - y$. Note that $J(C) $ acts as an abstract group on $\tilde A^\partial(\K) = \tilde A^\partial(\hat K)$, which we abbreviate as $\tilde A^\partial$. 
  
\begin{Theorem}Ê: let $x \in  L\tilde A(K)$, and let $B$ be the smallest abelian subvariety of $A$ such that   $\overline x$   lies in $LB+ LA_0(\C)$. Then, the Galois group  $Aut_\partial (\K(y)/\K)$ of $(**)$ over $\K$ is equal to $\tilde B^\partial$, and in particular, $tr.deg.(\K(y)/\K) = 2 dim B$. 
\end{Theorem} 
 
 \noindent
 {\it Proof} : by Theorem 3.2,  this Galois group is $H_a^\partial$, where $H_a$ is the smallest $D$-subgroup  $H/\K$  of $\tilde A$ such that
$a := \nabla_{L\tilde A} (x) \in LH + \partial \ell n_{\tilde A}(\tilde A (\K))$. Lifting $\overline x$ to a point $x'$ of $L\tilde B(K)$, we see by the same argument as above that $H_a \subset \tilde B$. To prove the reverse inclusion, we suppose for a contradiction that $G := \tilde B/H_a \neq 0$.  

\medskip
We first prove that $H_a$ is automatically defined over $K$. Here, we cannot appeal to the cohomological argument of Proof (i) of \S 3.2, because the full extension $\K(y)/K$ is not ``normal" in any reasonable sense, so that in the diagram above, there is no natural  action of $J(C)$ on $H_a^\partial$ by conjugation. Instead, we use  the rigidity of abelian varieties : the projection $A'' $ of $H_a$ on $A$ is necessarily defined over $K$, and $H_a$
is isogenous to $\tilde A'' \times W''$, where $W'' \subset W_A$ is a $D$-submodule of $L\tilde A$ over $\K$. Now,   $W''$  is generated over $\K$ by $W''^\partial$, which is stable under the action of $J(C)$ by the minimality of $H_a$. So, $W''$ is indeed defined over $K$, and $H_a$ as well.

\medskip
We now consider the $D$-quotient  $G = \tilde B/H_a$ over $K$, which, up to an isogeny, can also be viewed as a quotient of $\tilde A$, and denote by $\xi, \eta$ the images of $x, y$ in $ G$. Being stable under $Aut_\partial (\K(y)/\K)$, the point $\eta$ is  $\K$-rational, and the class of $\nabla_{ LG}(\xi) = \partial \ell n_{G}(\eta) \in LG(K)$ in $Coker(\partial \ell n_G , K)$ becomes trivial in $Coker(\partial \ell n_G, \K)$. Going to a minimal non trivial $D$-quotient $\overline G$ of $G$ over $K$, we will presently check that the natural map 
$$Coker(\partial \ell n_{\overline G} , K) \rightarrow Coker(\partial \ell n_{\overline G}, \K)$$
is injective. Consequently, $\eta$ lifts to a point $y' \in \tilde B(K)$ such that for some  proper $D$-subgroup $H'$ of $\tilde B$ over $K$,
$$\partial \ell n_{\tilde B}(y') \in  \nabla_{L\tilde B} (x') + LH' \subset \nabla_{L\tilde B}(L\tilde B(K)) + LH'.$$
From version ($1^*$) of Manin-Chai, we deduce that $\partial \ell n_{\tilde B}(y')$ lies in $\nabla_{L\tilde B}(W_B (K))+ LH'$, hence $ \nabla_{L\tilde B} (x')$ as well. In view of Facts 1  and 3, this contradicts the non-degeneracy of $\overline x$ in $LB$. 

\medskip
To prove the required injection of cokernels, we follow the principle of Proof (ii) of \S 3.2. But first, we note that since $H_a$ does not project onto $B$ by assumption,  the  quotient $\overline G$ is necessarily  a quotient of $\tilde B/U_B$ (in the notations of Fact 4), and is therefore $K$-large. So, we must in fact  show that  given two points $\overline \alpha \in L\overline G(K), \overline \eta \in \overline G(\K)$ such that $\overline \alpha = \partial \ell n_{\overline G}(\overline \eta)$, then $\overline \eta$ is automatically defined over $K$.  By Fact 2,  the abstract group $J(C) = Aut_\partial(\K/K)$ acts on $\tilde B^\partial = \tilde B^\partial (\K)$, inducing a trivial action on $\overline G^\partial(\K) =  \overline G^\partial(K)$ by $K$-largeness.  Therefore, the  cocycle 
$$ \hat \xi : J(C) \rightarrow \overline G^\partial(K) : \tau \mapsto \hat  \xi(\tau) := \tau(\overline \eta) - \overline  \eta$$
is a homomorphism of abstract groups. Since by Fact 3,  every element of $J(C)$ is a commutator, while $\overline G^\partial(K)$ is commutative, $\hat \xi$ is trivial.  In other words, $\overline \eta$ is fixed by $Aut_\partial(\K/K)$, and is therefore $K$-rational. (Another proof, suggested by A. Pillay, will be found in \cite {BM}, Remarque 7).  

\section{APPENDIX}  

Since $\tilde A$ disappears from the statements of this appendix,  we will now   call $y$ instead of $\overline y$  the points of $A$. We retain our  convention that a multiple of $y$ means a multiple by a non-zero integer. 

\subsection{Kummer theory}

Let $K$ be a {\it number field}, with algebraic closure $\overline K$, and let $A$ be an abelian variety over $K$. Going to a finite extension if necessary, we assume that all the endomorphisms of $A/\overline K$ are defined over $K$ and set
$$ End(A/K) = End(A/\overline K) := \cal O.$$
Let $y$ be a point in $A(K)$. Since there  is no ``constant part" anymore, we say that {\it $y$ is non-degenerate on $A$} if  for any proper abelian subvariety $B$ of $A$, no   multiple of $  y$  lies in $B$, i.e. if the group  ${\Z}.y$ generated by $y$ is Zariski dense  in $A$, or equivalently, if the annihilator $Ann_{\cal O}(y) $ of $y$ in $\cal O$ is reduced to $\{0\}$.  

\medskip
Following the elliptic work of Bashmakov and Tate-Coates from the 70's (see \cite{L}, V, \S 5),  K. Ribet devised in \cite{Ri} a general method to bound from below the degree of the division points of $y$ of prime order. His method readily extends to all orders (see \cite{BK}), and yields :

\begin{Theorem} :  let $y$ be a non-degenerate point in $A(K)$. There exists a real number $c = c(A,K,y) > 0$ such that for all $n > 0$, $$[K({1 \over n} y) : K] \geq c n^{2 dimA}.$$
\end{Theorem}
Set $dim A := g$, and write $[n]\in {\cal O}$ for the multiplication by $n$ on $A$, with kernel $A[n] \simeq (\Z/n\Z)^{2g}$ in $A(\overline K)$. 
Since $A$ is usually not ``$K$-large" for $[n]$,   we introduce a field $K_\infty$, analogous to the previously defined field $\K$, which takes into account all torsion points of $A$ :
$$K_\infty= K(A_{tor}), ~{\rm~where~} A_{tor} = \cup_{n>0}ÊA[n].$$
For each  $\ell$ in the set $\cal P$ of prime numbers, we set
$$    K_{y,(\ell)} = \cup_m K_\infty({1 \over \ell^m} y) , \quad K_{y,\infty} = \cup_{n} K_\infty({1 \over n} y), $$
and define the Tate modules $  T_\infty(A) := \varprojlim_{~n} A[n]  = \Pi_{\ell \in {\cal P} } T_\ell(A)$ in the usual way. 

\medskip
Given a positive integer $n$, we have the tower of Galois extensions of $K$, whose Galois groups are indicated in the diagram :

  \begin{equation*}
\begin{array}{cccccc}
&&K_\infty({1\over n} y) & &\xi_y \\
&& \mid &  \} N  & \hookrightarrow & A[n] \simeq ({\Z}/n{\Z})^{2g}\\
&\Gamma \{ &K_\infty & &  \rho \\
&& \mid & \} J & \hookrightarrow & GL(T_\infty(A)) \simeq GL_{2g}(\hat \Z) \\
&& K  \\
 \end{array}
 \end{equation*}
 
\noindent
where $N \ni \sigma \mapsto \xi_y(\sigma) = \sigma({1\over n} y) - {1 \over n} y$ is a group embedding. Theorem 5.1 immediately follows from :
\begin{Theorem} :  let $y$ be a $K$-rational point of $A$, and let $B$ be the smallest abelian subvariety of $A$ containing a multiple  of $y$. Then, $Gal(K_{y,\infty}/K_\infty)$ is isomorphic to  an open subgroup of $T_\infty(B)$.
\end{Theorem}

\noindent
Since all torsion points are defined over $K_\infty$, we can  assume that $y$ itself lies in $B$; a map analogous to $\xi_y$ then identifies this Galois group to a subgroup of  $T_\infty(B)$.  By  Nakayama, Theorem 5.2 then amounts to   the following claims :

\medskip

i) for  almost all $\ell \in {\cal P}$, $Gal(K_\infty({1\over \ell} y)/K_\infty) \simeq B[\ell] \simeq (\Z /\ell \Z)^{2 dimB}.$

ii) for  all $\ell \in{\cal P}$, $Gal(K_{y, (\ell)}/K_\infty)$ is an open subgroup of $T_\ell(B) \simeq  {\Z}_\ell^{2dimB}$;

\subsection{\bf Ribet's method.}  

The proof  ``follows" that of Theorems 3.3, 4.3, 4.4. For the sake of simplicity, we present it under the assumption that  $y$ is non-degenerate, i.e. $B = A$, and refer to \cite{Ri} for the general case  (see also \cite{H}, \cite{BK}).  We fix a prime number $\ell$, and   first treat the mod $\ell$ claim.

\medskip
\noindent
1. {\it  Galois theoretic step} ($=$ analogue of ``Proof (i)" of \S 3.2). 

\medskip
By the same computation from affine geometry, $Im(\xi_y) \simeq N$ is a $J$-submodule of $A[\ell]$. In these conditons, there exists $\ell_1 = \ell_1(A,K)$ with the following property : assume that $N \neq A[\ell]$ and that  $\ell > \ell_1$; then, there exists $\alpha \in {\cal O}$, $ \alpha \notin \ell {\cal O}$, such that  $\alpha. y$ is divisible by $\ell$ in $A(K_\infty)$.

\medskip
\noindent
Indeed, Faltings's theorem on the finiteness of isomorphism classes of abelian varieties over $K$ in a given isogeny class implies that for $\ell$ large enough, $N$ is the intersection with $A[\ell]$ of the kernel of an endomorphism  $\alpha$ of $A$. Alternatively, one can say, as in \cite{Ri}, that :

- $A[\ell]$ is a semi-simple $J$-module, so there exists $\alpha_\ell \in End_J(A[\ell]), \alpha_\ell \neq 0$, killing $N$.

- $End_J(A[\ell]) \simeq End(A) \otimes {\bf F}_\ell$ (Tate conjecture), so $\alpha_\ell$ is induced by some  $\alpha \in {\cal O}, \alpha \notin \ell {\cal O}$. 

\noindent
Finally,

- $\xi_{\alpha.y} = \alpha \xi_y$, so, ${1 \over \ell} \alpha. y $ is fixed by $N$, i.e. the class of $\alpha. y$ in $Coker([\ell], A(K_\infty))$ vanishes, where we set, for any extension $K'$ of $K$:
$$Coker([\ell], A(K')) := A(K') /\ell A(K').$$

\medskip
\noindent
{\bf Remark 1} :  in the function field case, and even if $A$ is simple and non isoconstant, the hypothesis that $(End(L\tilde A))^{\nabla_{L\tilde A}} = End(A) \otimes \C$, which  would be  a functional  version of Tate's conjecture, {\it does not hold} in general. See \cite{De} \S 4.4,  and a counterexample in \cite{Fa}. 

\medskip
\noindent
2. {\it  Galois descent} ($=$ analogue of ``Proof (ii)" of \S 3.2)

\medskip
There exists $\ell_2(A,K)$ such that if $\ell > \ell_2$, and  if a point $y'  \in A(K)$ is divisible by $\ell$ in $A(K_\infty)$, then, $y'$ is already divisible by $\ell$ in $A(K)$, i.e. the kernel $Ker$ of the natural map
$$Coker ([\ell], A(K) ) \rightarrow Coker ([\ell], A(K_\infty))$$
vanishes.

\medskip
\noindent
Indeed, the exact sequence of inflation-restriction here gives:

  \begin{equation*}
\begin{array}{cccccc}
Ker &  \rightarrow &  A(K) /\ell. A(K)   & \rightarrow  &A(K_\infty)/\ell. A(K_\infty)\\
  \downarrow & &\downarrow {\hat \xi}  & & \downarrow {\xi}\\
H^1(J, A[\ell])& \rightarrow &  H^1(\Gamma, A[\ell]) & \rightarrow & Hom_J(N, A[\ell]) ~, \\
\end{array}
\end{equation*} 
while  Serre's results on homotheties \cite{S}, Thm. 2, and Sah's lemma imply that $H^1(J, A[\ell]) = 0$ for $\ell$ large enough (depending only $A$ and $K$).  

\medskip
\noindent
3. {\it  Diophantine   step}

\medskip
By Step 2,  $\alpha.y$ is divisible by $\ell$ in $A(K)$.
A contradiction to the conclusion of Step 1,  or to the non-degeneracy of $y$, now follows from the existence of $\ell_0 = \ell_0(A, K, y)$ such that if $\ell > \ell_0$, and if $\alpha. y \in \ell.A(K)$, then there exists $\alpha' \in {\cal O}$ such that $(\alpha - \ell \alpha').y = 0$.   

\medskip
\noindent
Indeed, this in turn  follows from the  Mordell-Weil theorem, which implies that ${\cal O}.y$ has finite  index in its divisible hull   in $A(K)$  (see \cite{L} for an effective version). Since Manin's kernel theorem is  based on a similar fact,  this last step can perhaps be considered  as an analogue of Formulae (1), (2), ($1^*$), ($2^*$) of the text. 

\medskip
\noindent
{\bf Remark 2} : even when $y$ is ``indivisible in $A(K)$" in the sense of  \cite{H}, Lemme I,  the constant $\ell_0(A, K, y)$ arising in this last step {\it cannot be made independent of $y$} in general. Elliptic curves $A$ with complex multiplications by non-principal orders $\cal O$ already provide counterexamples. Consequently, the constant $c$ occurring in Theorem 5.1 will in general depend on  the non-degenerate point $y$, even if one insists that ${\cal O}.y$ be maximal among the $\cal O$-orbits of points of $A(K)$. 

\bigskip

The $\ell$-adic claim can be treated along similar lines\footnote{~See \cite{BK}, Theorem 2, for a very brief sketch. As pointed out to me by M. Bays, the argument is described in more detail in \cite{BGK}.}, as follows. Firstly, by \cite{S}, Thm. 1, we may assume, after a finite extension of $K$, that  the fields of definition $K_{(p)}, p \in {\cal P},$ of the $p$-primary parts of $A_{tor}$ are linearly disjoint over $K$. Since $\Z_\ell$ cannot be a quotient of a  $p$-adic Lie group for $p \neq \ell$, it  then suffices to prove Claim (ii) with $K_\infty$ replaced by $K_{(\ell)}$ and $K_{y,(\ell)} $ by $K_{(\ell), y} := \cup_m K_{(\ell)}({1 \over \ell^m} y)$. So, we have a continuous map $\xi_{(\ell),y}$ analogous to $\xi_y$, and must show that $\xi_{(\ell),y}$ sends $N_{(\ell)} := Gal(K_{(\ell),y}/K_{(\ell)})$  into an open subgroup of $T_\ell(A)$. Set $J_{(\ell)} = Gal(K_{(\ell)}/K)$, and denote by $\Xi_{\ell} : A(K) \otimes \Z_\ell \rightarrow Hom(N_{(\ell)}, T_\ell(A))$  the $\Z_\ell$-linear extension of the map $y \mapsto \Xi_\ell(y) = \xi_{(\ell),y}$. 

\medskip
\noindent
1. {\it  Galois theoretic step}  

\medskip
 $N_{(\ell)}$ is as usual a $J_{(\ell)}$-submodule of $T_\ell(A)$,  and is closed. Assuming for a contradiction that it is not open, we deduce from the semi-simplicity of the representation $T_\ell(A) \otimes \Q_\ell$  and from Tate's conjecture the existence of a non-zero element $\alpha \in {\cal O}\otimes \Z_\ell$,  such that  $\Xi_\ell(\alpha. y) = 0$.

\medskip
\noindent
2. {\it  Galois descent} 

\medskip
By the inflation-restriction sequence, the kernel of $\Xi_\ell$ injects into $H^1(J_{(\ell)}, T_\ell (A))$. Bogomolov's theorem on homotheties \cite{Bog} or, more directly, an earlier result of Serre on the vanishing of $H^1(J_{(\ell)}, T_\ell (A)\otimes \Q_\ell)$  ensure that the latter group is finite. Replacing $\alpha$ by some multiple, we deduce that $\alpha.y = 0$ in $A(K) \otimes \Z_\ell$.
 
\medskip
\noindent
3. {\it Diophantine   step}

\medskip
Choosing a basis of $\cal O$ over $\Z$, and a basis of $A(K)$ over $\Z$ modulo torsion, we deduce from the latter conclusion that $y$ is linearly dependent over $\cal O$, contrary to our assumption that $y$ is non-degenerate in $A(K)$.

\subsection{Back to function fields}
 Let us come back to the situation  of  \S 4, with a base field $K = \C(S)$ and an abelian variety $A$ over $K$, with $\C$-trace $A_0$. All the notions introduced in the present Appendix remain meaningful, and one can ask if a suitable version of Theorem 5.2 still holds true in this functional context.
 
 This is indeed the case. Although this probably follows from  Theorem 5.2 itself and a specialization argument, we here want to indicate a more natural method, where the searched-for  algebraic  statement on {\it classical} Galois groups (of finite extensions of $\C(S)$) is   {\it  directly deduced}  from the purely transcendental Theorem 4.3. on {\it differential} Galois groups (of Picard-Vessiot extensions of $\C(S)$).
 
\medskip
Since $\tilde A$ will reappear, we start with a point $\overline y$ in $A(K)$,   denote by $B$ the smallest abelian variety such that a non trivial  multiple of $\overline y$ lies in $B + A_0(\C)$,  assume without loss of generality that $\overline y$ itself lies in $B$, and choose an arbitrary lift $y$ of $\overline y$ to $\tilde B(K)$. Fix a point $s_0 \in S(\C)$,  and consider the image $\Pi$ of  $\pi_1(S, s_0)$ in the monodromy representation attached to the local system  formed by the various ``logarithms" of   the multiples $my, m \in \Z$, of $y \in {\cal \tilde B}(S)$. More precisely, let $x = \ell n_{\tilde B}(y) \in L\tilde B$ be a local determination of a logarithm of $y$ in a neighbourhood of $s_0$. For any $\gamma \in \pi_1(S, s_0)$, analytic continuation along $\gamma$ provides an element of the differential Galois group 
$$\Gamma(\C) = Aut_\partial(\K(x)/K), {\rm ~Êwith}~\Gamma \subset GL_{2g+1}, ~{\rm and}~\Pi \subset \Gamma(\Z) \subset GL_{2g+1}(\Z),$$ sending $x$ to $\gamma.x = x + \hat \xi(\gamma)$, where $\hat \xi(\gamma)$ lies in the subgroup $\Omega_{\tilde {\cal B}}$ of periods  of $(L\tilde B)^\partial$. Notice that $\hat \xi$ is only a	 cocyle, so that $\hat \xi(\pi_1(S, s_0)) := \Omega_y$ is in general not a group.

\medskip
Let now $n$ be a positive integer, and consider the $n$-th division point ${1 \over n} y =  exp_{\tilde {A}}({1 \over n} x) $ of $y$ in $\tilde B(\overline K)$.   Since $exp_{\tilde { A}}$ is a $S^{an}$-morphism, the action of $\pi_1(S, s_0)$ on ${1 \over n} y$  is given by 
$$\gamma.({1 \over n} y) = exp_{\tilde { A}}(\gamma.({1 \over n} x)) =   exp_{\tilde A}({1\over n} x + {1 \over n} \hat \xi(\gamma)) = {1 \over n} y +  exp_{\tilde {B}}( {1 \over n} \omega_\gamma ),$$
 where $\omega_\gamma = \hat \xi(\gamma) \in \Omega_y$.  In particular the number of conjugates of ${1 \over n} y$ over $K$, i.e. the degree of $K({1 \over n} y)$ over $K$,  is equal to the number $\delta_y(n)$ of distinct classes
 in  ${1 \over n}Ê\Omega_y$ modulo the kernel $ \Omega_{\tilde {\cal B}}$ of $exp_{\tilde { B}}$. And since $\tilde B$ is a vectorial extension, the fields of definition   of ${1 \over n}y$ and of its image ${1 \over n}\overline y$ in $B$ coincide, so this is also the degree of $K({1 \over n} \overline y)$ over $K$.

\bigskip

Now come the main points : 

\medskip

(i) since $\nabla_{L\tilde A}$ is a fuchsian connexion, $\Pi$ is a Zariski-dense subgroup of the algebraic group $\Gamma$. In particular, $\Omega_y$ is Zariski dense in the Galois group $(N_a^\partial) $, which, as we know by Theorem 4.3,  coincides with $(L\tilde B)^\partial$. 

\medskip

(ii) Assume for a moment that $A$ is defined over $\C$, i.e. that $A = A_0$, in which case all the torsion points of $A$ are defined over $\C$, and $K_\infty = K $. Then, the periods in $ \Omega_{\tilde {\cal A}}$ are constant (so, $K = \K$ as well), $\hat \xi = \xi$ is a group morphism, and $\Omega_y$ is a subgroup of the discrete subgroup   $\Omega_{\tilde {\cal B}}$ of the vectorial group $(L\tilde B)^\partial$. Being Zariski-dense in the latter vectorial group, $\Omega_y$ must be of {\it finite index}, say $\nu_y$, in $\Omega_{\tilde {\cal B}}$. In particular, the degree $\delta_y(n) = [K({1 \over n} \overline y):K] $ is bounded from below by $c n^{2 dim B}$, where $c = {1 \over \nu_y}$, and  we derive, more precisely : 
 
\begin{Theorem} let $A_0$ be an abelian variety defined over $\C$, let $K = \C(S)$,  let $y$ be a point in $  A_0(K)$, and let $B$ be the smallest abelian subvariety of $A_0$ such that  some multiple of  $  y$   lies in $B+ A_0(\C)$. There exists an integer  $\nu = \nu(A_0,K, y)$ such that for any $n > 0$,  the Galois group  $Gal (K({1\over n}  y))/K)$ is isomorphic to a subgroup of $B[n]$ of bounded index, equal to 1 as soon as $n$ is prime to $\nu$.  
\end{Theorem}

(iii) When $A$ is not defined over $\C$, the discrete monodromy group $\Pi$ will in general not even meet $N_a^\partial$ outside of $0$, and the above argument does not apply. However, finitely generated subgroups of $GL_m(\Z)$ such as $\Pi$ often satisfy the {\it strong approximation  property} with respect to their Zariski closure $  G$ in $GL_m$, in the sense   that their closure $\hat \Pi$ in the profinite group $GL_n(\hat \Z)$ is then open in $G(\hat \Z)$.   This holds true when $G$ is a semi-simple connected and simply connected group, as shown by the theorems of Matthews-Vaserstein-Weisfeiler and of Nori, which play a role in \cite{S}. It clearly fails for tori, but as pointed out to me  by Y. Benoist,  Nori's Theorem 5.3 in \cite{No} shows that this is in a sense the only obstruction\;: by Fact 3 of \S 4, this theorem applies to our group $G = \Gamma$, whose radical $N_a^\partial$ is unipotent. We deduce that for all prime numbers $\ell$, the image of $\Pi$ in $\Gamma(\Z/\ell \Z)$ has bounded index. In particular\footnote{~Actually, this consequence follows directly from Nori's Theorem 5.1 in \cite{No}, according to which for almost all $\ell$'s, the image of $\Pi$ in  the group $\Gamma(\F_\ell)$  contains the subgroup  $\Gamma(\F_\ell)^+$ generated by its elements of order $\ell$, hence the full subgroup  $N_a^\partial(\F_\ell) \simeq \F_\ell^{\,2 dim B}.$}, the image  of 
 ${1 \over \ell}Ê\Omega_y$ in   ${1 \over \ell} \Omega_{\tilde {\cal B}}/ \Omega_{\tilde {\cal B}} \simeq B[\ell]$ fills up this group for $\ell$ sufficiently large. 
 
 \medskip
 We will come back to this in a later article, but already mention the following  consequence of the discussion above, where we set  $K_n = K(A[n])$ :
  
\begin{Theorem}Ê: let $A$ be an abelian variety over $K = \C(S)$, with $\C$-trace $A_0$, let $y$ be a point in $  A(K)$, and let $B$ be the smallest abelian subvariety of $A$ such that  some multiple of  $  y$   lies in $B+ A_0(\C)$. There exists an integer $\nu =  \nu(A,K, y)$   such that for any $n > 0$,  the Galois group  $Gal (K_n({1\over n}  y))/K_n)$ contains a subgroup of $B[n]$ of bounded index, and coincides with $B[n]$  as soon as $n$ is prime to $\nu$.  

\end{Theorem}

\bigskip
\centerline{***}

\bigskip
\noindent
{\it Adresse de l'auteur} : bertrand@math.jussieu.fr

\medskip
\noindent
{\it Mots clefs }: linear and non-linear differential Galois theory; abelian varieties; Galois cohomology;  Kummer theory.

\medskip

\noindent
{\it  Classification AMS} : 12H05, 14K15, 11G10, 12G05. 
\end{document}